\title{Continuity of the SLE trace \\ in simply connected domains}
\author{Christophe Garban\footnote{Microsoft Research, Universit\'e Paris-Sud and ENS, partially supported by the
ANR under the grant ANR-06-BLAN-0058.}
\and Steffen Rohde\footnote{University of Washington, Research supported in part by NSF Grants DMS-0501726
and DMS-0800968.}
\and Oded Schramm\footnote{Microsoft Research.}}

\documentclass[12pt]{article}

\usepackage{amsfonts,amsthm,amsmath,latexsym,amssymb,graphicx}

\input epsf.sty

\newcommand{\diam}{\operatorname{diam}}

\newtheorem{thm}{Theorem}[section]
\newtheorem{prop}[thm]{Proposition}

\newtheorem{lemma}[thm]{Lemma}

\def\D{{\mathbb D}}
\def\R{{\mathbb R}}
\def\H{{\mathbb H}}
\def\C{{\mathbb C}}
\def\E{{\mathbb E}}
\def\P{{\mathbb P}}
\def\LL{{\cal L}}
\def\QQ{{\cal Q}}
\def\eps{\varepsilon}

\def\1{{\bf 1}}
\def\p{\partial}
\renewcommand{\dim}{\operatorname{dim}}

\begin{document}

%\baselineskip=1.2\baselineskip
%\date{\today}
\date{}
\maketitle

\begin{abstract}We prove that the $\mathrm{SLE}_\kappa$ trace in any simply connected domain
$G$ is continuous (except possibly near its endpoints)
if $\kappa<8$. We also prove an SLE analog of Makarov's Theorem about the support of harmonic measure.
\end{abstract}

%\begin{doublespace}
\section{Introduction}

The stochastic Loewner evolution (SLE) describes a collection of random curves that are
related to scaling limits of two-dimensional statistical physics systems.
In~\cite[Theorem 5.1]{RS} it was shown that
the chordal SLE trace in the upper half plane $\H$
is a well defined continuous path. For other simply connected domains
$G\subsetneqq\C$, the SLE in $G$ is defined via a conformal homeomorphism
$f:\H \to G$. The situation with radial SLE is similar, except
that the \lq\lq standard\rq\rq\ domain is the unit disk $\D$.
Our first theorem extends this continuity result, as follows.

\begin{thm}\label{continuity}
Let $G\subsetneqq \C$ be a simply connected domain, let $a,b$ be two prime ends
of $G$, let $z_0\in G$, and let $\kappa\in[0,8)$.
Then the chordal $\mathrm{SLE}_\kappa$ trace in $G$ from $a$ to $b$
and the radial $\mathrm{SLE}_\kappa$ trace in $G$ from $a$ to $z_0$
are a.s.\ continuous in $(0,\infty)$.
\end{thm}

Besides the intrinsic interest in this result, it is also useful in the
general theory of SLE and the related scaling limits. For example, 
the construction of the conformal loop ensembles in~\cite[Section 4.1]{SheffCLE}
would have been simpler if this theorem was available.

If the boundary $\partial G$ is a smooth curve
(more generally, if it is locally connected), then the conformal
map $f$ to $G$ extends continuously to the closure of $\H$ (respectively $\D$),
and the continuity of the trace $f\circ\gamma$ follows at once from the continuity 
of the trace $\gamma$. But if $\partial G$ contains boundary points at
which it looks like the topologist's sine curve, then $f$ is not continuous at the corresponding
points, and the continuity of  $f\circ\gamma$ is no longer obvious when $\kappa>4$.
In fact, this non-continuity could happen at {\it every} boundary point: there are simply
connected domains $G$ for which the limit set 
$f(z):=\{w: \exists (z_k\to z), \lim f(z_k)=w\}$ equals $\partial G,$ for {\bf all} $z\in\partial\D$ 
\cite{K}!  On the other hand, for every conformal map $f:\D\to\C$, 
the radial limit $lim_{r\to1} f(r e^{it})$ exists for a.e.\ $t\in[0,2\pi]$. A celebrated theorem
of Makarov asserts that there is a set $A\subset\partial\D$ of full measure such that the set
of radial limits $f(A)$ has Hausdorff dimension 1 (even sigma-finite length). 
Equivalently, for every simply connected domain $G\subsetneqq \C$ there is a set $B\subset\partial G$ 
of Hausdorff dimension 1 such that a Brownian motion started inside $G$ will a.s.\ exit $G$ through $B$.
However, under a mild assumption on the geometry of $G$ (precisely, if $G$ is a John domain),
reflected Brownian motion in $G$ intersects the boundary in a set of full dimension \cite{BCR}. 
In particular, there is no nontrivial upper bound on the dimension of the trace of reflected Brownian motion 
on $\partial G$. The situation is different for SLE:

%%%%%%%%%%%%%%%%%%
\begin{figure}[htbp]
\centerline{\epsfysize=3 true in \epsffile{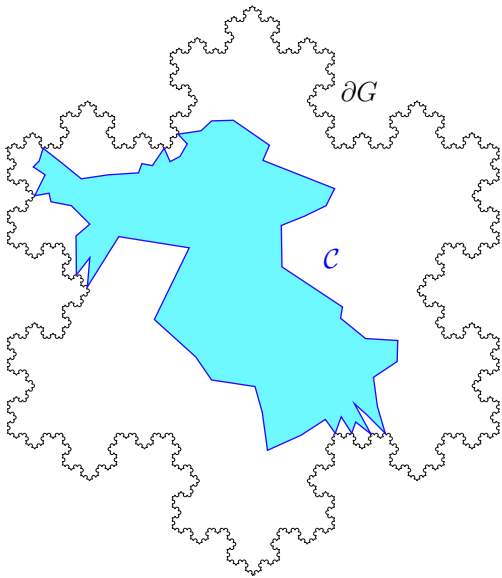}} 
\caption{\label{f.percocase} A schematic view of a percolation cluster $\mathcal{C}$
(or an $\mathrm{SLE}_6$ hull) inside a fractal domain $G$; the blue curve represents the exterior boundary of the cluster.}
\end{figure}
%%%%%%%%%%%%%%%%

\begin{thm}\label{dimension} 
Let $G\subsetneqq \C$ be a simply connected domain, let $a,b$ be two prime ends
of $G$, let $z_0\in G$, and let $\kappa\in(4,8)$.
Then there is a Borel set (actually a $F_\sigma$ set) $B\subset\partial G$ such that 
the chordal $\mathrm{SLE}_\kappa$ trace in $G$ from $a$ to $b$
and the radial $\mathrm{SLE}_\kappa$ trace in $G$ from $a$ to $z_0$
almost surely satisfy
$$\gamma(0,\infty)\cap\partial G \subset B,$$
and 
$$\dim B \leq d(\kappa)<2 \,,$$
where $d(\kappa)$ is a constant that depends only on $\kappa$.
\end{thm}

In the case $\kappa=6$, this theorem can be thought of as a Makarov theorem for percolation.
Indeed, instead of starting a Brownian motion at $z_0$ inside $G$, think of ``sending'' a critical percolation 
cluster from $z_0$ in the following way: already at the scaling limit, condition $z_0$ to be connected via some 
open cluster to the boundary $\p G$ (though this is an event of probability 0, it is possible to make sense of 
this conditioning; see \cite{Ke}).
Theorem \ref{dimension} then implies that there is a 
subset $B\subset\p G$ of dimension $d\leq d(6)<2$ which almost surely
``absorbs'' all the points on $\p G$ which are connected to $z_0$ within $G$. See figure \ref{f.percocase}
for an illustration of this. 
Of course one cannot hope to find such a set $B$ of dimension one, as in Makarov's theorem, since the 
percolation cluster is much ``thicker'' than the Brownian path stopped on the boundary.
Following this intuition, we will also show that the smallest bound $d(\kappa)$ necessarily satisfies $d(\kappa)>1$ for all
$\kappa \in (4,8)$.

The case $\kappa=4$ should nearly correspond to the setting of Makarov's theorem (for instance if one considers chordal $\mathrm{SLE}_4$ in the disc with 
a random initial point). Therefore we expect 
that for $\kappa$ close to 4, $d(\kappa)$ should be close to one. We indeed prove the following estimates on
$d(\kappa)$:

\begin{prop}\label{p.estimate}
There are absolute constants $C_1,C_2>0$ such that for any $\kappa\in (4,8)$,
\begin{equation}
d(\kappa) \leq  \bigl(2 - C_1 (8-\kappa)\bigr) \wedge \bigl(1+  C_2 \sqrt{ \kappa -4 }\bigr)\,.
\end{equation}
\end{prop}

Finally, in Section 4 we will relate the integral means spectrum of a conformal map to
the dimension of the SLE trace on the boundary of a domain, and show existence of nice
Jordan curves such that two independent $\mathrm{SLE}_\kappa$, run in the two complimentary
domains of the curve, are almost surely disjoint.\\

\noindent{\bf Acknowledgments}.
We wish to thank Jeff Steif for pointing out an error in the statement of Theorem \ref{dimension}.

\section{Uniform Continuity}

For sets $A\subset \C,$ denote
$$H_p(A)= \inf\Bigl\{ \sum r_i^p  : A\subset \bigcup_i B(x_i, r_i)\Bigr\},$$
which is called
the $p$-dimensional Hausdorff content  of $A$, where the infimum is
over all covers of $A$ by discs with positive radius. The following is an adaptation 
of \cite[Proposition 3.3]{KR}.

\begin{lemma}\label{content}
Let $G\subsetneqq\C$ be a simply connected domain and
$f:\D\to G$ a conformal homeomorphism. For every $0<p<1$ and $\eps>0$
there is $D\subset\D$ and $C>0$ such that
$$H_p(\D\setminus D) <\eps$$
and 
\begin{equation}
\label{e.ho}
|f(z)-f(z')| \leq C\, |z-z'|^{p/2}
\end{equation}
for all $z,z'\in D.$
\end{lemma}

For our present application, the case where $p$ is small is
the most relevant. The proof shows that we can choose the centers
of the discs to lie on the unit circle.

\begin{proof}
We first assume that $G$ is bounded.
Consider the collection $\QQ$ of dyadic \lq\lq squares\rq\rq\ 
$$Q=Q_{n,k} = \Bigl\{ r e^{i t} : 1-2^{-n}\leq r < 1, \frac k{2^n}\leq \frac t{2\pi} \leq \frac {k+1}{2^n}\Bigr\},$$ 
where $n\geq1$ and $k=0,1,2,...,2^n-1.$
Denote by $\ell(Q):=2^{-n}$, the size of $Q$,
and $T(Q):=\{z\in Q: 1-2^{-n}\leq r \leq 1-2^{-(n+1)}\}$,
the \lq\lq inner half\rq\rq\ of $Q$.
Fix $N>1$, to be determined later, and let $\LL$ be the collection of those $Q\in\QQ$ for which
$$\ell(Q) \leq 2^{-N}$$ and
$$\int_{T(Q)} |f'|^2 > \ell(Q)^p.$$ 
Set
$$D=\D\setminus \bigcup_{Q\in \LL} Q\,.$$ 
We claim that
\begin{equation}\label{derivative}
\forall_{z\in D}\qquad
|f'(z)| \leq C\, \frac1{(1-|z|)^{1-p/2}}
\end{equation}
with $C$ depending on $N$ and $p$ only.
Indeed, let $z\in D$. If $|z|<1-2^{-N}$,
we get~(\ref{derivative}) simply by choosing
$C$ large enough. Else, suppose that $|z|\geq1-2^{-N}$.
Let $Q$ be such that $z\in T(Q)$ and notice that $Q\notin\LL$. 
Hence $\int_{T(Q)} |f'|^2 \leq \ell(Q)^p.$
By the Koebe distortion theorem~\cite{P}, $|f'|$ is essentially constant in $T(Q)$
and (\ref{derivative}) follows.

We now claim that $f$ satisfies~\eqref{e.ho} for $z,z'\in D$, with possibly a 
different constant $C$.
Consider any $z,z'\in D$. First, suppose that $z=s\,z'$, where $s>1$.
Note that the interval $[z,z']$ is contained in $D$. Then we may
integrate the estimate~\eqref{derivative} over $[z,z']$ to
obtain~\eqref{e.ho}. 
Next, suppose that $z=r\,e^{i\theta_1}$ and $z'=r\,e^{i\theta_2}$,
where $|\theta_1-\theta_2|/(2\,\pi)\le 1-r$.
In that case, the path $r\,e^{i\theta}$, $\theta\in[\theta_1,\theta_2]$,
is contained in the union of some $Q\in\QQ$ satisfying $z\in T(Q)$
and a possibly different $Q'\in\QQ$ satisfying $z'\in T(Q')$. Therefore, this path is in $D$,
and we get~\eqref{e.ho} in the same way.
In general, suppose that $z=r_1\,e^{i\theta_1}$ and $z'=r_2\,e^{i\theta_2}$
with $|\theta_1-\theta_2|\le \pi$.
Then take $\rho:=\min\{r_1,r_2,1-|\theta_1-\theta_2|/(2\,\pi)\}$.
We then use the above cases and
$|f(z)-f(z')|\le |f(z)-f(\rho\,e^{i\theta_1})|+
|f(\rho\,e^{i\theta_1})- f(\rho\,e^{i\theta_2})|+
|f(\rho\,e^{i\theta_2})- f(z')|$, to obtain~\eqref{e.ho}.

To estimate the $p$-content of $\D\setminus D,$ just notice that the
interiors of the sets $T(Q)$, $Q\in\QQ$, are disjoint and
$$\sum_{Q\in\LL} \ell(Q)^p \leq \sum_{Q\in\LL} \int_{T(Q)} |f'|^2
\leq \operatorname{area}\{ f(z):1-2^{-N} < |z| <1\}\,,$$
which can be made arbitrarily small by choosing $N$ large.
This completes the proof in the case where $G$ is bounded.

The case of unbounded $G$ requires a few minor adaptations.
Set
$$
\phi(z):= \Bigl|\frac{f'(z)}{\max\bigl\{|f(z)|\log|f(z)|,1\bigr\}}\Bigr|^2,
$$
and redefine $\LL$ to be the set of $Q\in\QQ$ such that $\ell(Q)\le 2^{-N}$ and
$\int_{T(Q)} \phi>\ell(Q)^p$. Note that a simple change of variables gives
$$
\int_{\D}\phi=\int_G \max\{|z|\log |z|,1\}^{-2}<\infty
$$
and so we get $\sum_{Q\in\LL} \ell(Q)^p<\eps$ by taking $N$ sufficiently large,
as above.  The Koebe distortion theorem implies that
$\phi$ is essentially constant in $T(Q)$. Therefore,
we get 
\begin{equation}
\label{e.newbd}
\forall _{z\in D}\qquad
\Bigl|\frac{f'(z)}{\max\bigl\{|f(z)|\log|f(z)|,1\bigr\}}\Bigr|\le C\, \bigl(1-|z|\bigr)^{-1+p/2}.
\end{equation}
Set
$$
g(z):=\int_0^{|z|} \frac{ds}{\max\{s\log s,1\}}\,.
$$
Then $\frac d{dr}\, g\bigl(f(r\,e^{i\theta})\bigr)$ is bounded by the left hand side
of~\eqref{e.newbd} with $z=r\,e^{i\theta}$. Therefore, since $D$ is star-shaped about
$0$, we have for $z\in D$,
$$
g\circ f(z)\le g\circ f(0)+ C\int_0^1 \bigl(1-r\bigr)^{-1+p/2}\,dr<\infty\,.
$$
 Thus $g\circ f$ is bounded on $D$. Since $\lim_{z\to\infty}g(z)=\infty$,
it follows that $f$ is bounded on $D$.
Thus~\eqref{e.newbd} implies~\eqref{derivative} with possibly
a different constant. This gives~\eqref{e.ho}, as before.
\end{proof}

%\end{doublespace}

\section{Remaining proofs}

\noindent
{\bf Proof of Theorems \ref{continuity} and \ref{dimension}.}
Let $\gamma$ be the chordal $\mathrm{SLE}_{\kappa}$ trace in $\D$ from $1$ to $-1,$ 
and denote by $f$ a conformal map from $\D$ to $G$ sending $1$ and $-1$ to $a$ and $b$.
If $\kappa\leq4$, then $\gamma$ is continuous on $(0,\infty)$ and 
$\gamma(0,\infty)\subset \D$ almost surely  \cite[Theorems 5.1 and 6.1]{RS}.
The continuity of $f\circ\gamma$ follows at once. Now let $4<\kappa<8$, let
$\delta,\eps>0$ and $\delta<|t|<\pi-\delta$. 
We have
\begin{equation}
\label{hitting}
\P\bigl[\gamma(0,\infty)\cap B(e^{it},r)\neq\emptyset\bigr] \leq C\, r^{\frac8{\kappa} - 1}
\end{equation}
for some constant $C=C(\kappa,\delta)$ and for all $r<\delta/2$; see, e.g.,~\cite[Proposition 2.3]{SZ}
or~\cite[Theorem 3.2]{AK}.
Let $p=\frac8{\kappa} - 1$ and let $D$ be as in Lemma \ref{content}.
It follows from that lemma that there are discs $B(x_i,r_i)$ with $\D\setminus D\subset \bigcup_i B(x_i,r_i)$
and $\sum r_i^p<\eps,$ and we may and will assume $x_i\in\partial\D$.
From (\ref{hitting}) we obtain
$$\P\Bigl[\bigl(\gamma(0,\infty)\setminus (B(1,2\,\delta)\cup B(-1,2\,\delta)) \bigr)\cap \bigcup_i B(x_i,r_i))\neq\emptyset\Bigr] \leq C \eps.$$
Since $f$ is continuous on $\overline D$, it follows that with probability at least $1-C\eps$
$f\circ \gamma$ is continuous at every $t$ such that $\gamma(t)\notin B(1,2\,\delta)\cup B(-1,2\,\delta)$.
Now Theorem \ref{continuity} in the chordal case follows by first letting $\eps\to 0$, then 
letting $\delta\to 0$, and using the 
transience of $\gamma$~\cite[Theorem 7.1]{RS}.
The radial case is similar.

\medskip

To prove Theorem \ref{dimension}, consider the domain 
$$D_{\eps}=\D\setminus \bigcup_i T(x_i,r_i),$$ 
where $x_i$ and $r_i$ are as above and $T(x,r)$ is the triangular
region bounded by the circular arc $B(x,2\,r)\cap\partial\D$ and the two line segments joining
the endpoints of the arc with the point $(1-2\,r)\,x$. (Alternatively, let $D_{\eps}$ be the domain
$D_{\eps}=\D\setminus \bigcup_{Q\in \LL} Q$ from Lemma \ref{content}.) 
Then $D_{\eps}$ is a John domain (meaning that $D_\eps$ is bounded and
there is a number $M>1$ such that every Jordan 
arc $\gamma\subset D_{\eps}$ with endpoints on $\partial D_{\eps}$
decomposes $D_{\eps}$ into two subdomains at least one of which has diameter $\leq M \diam \gamma$)
with uniformly bounded John constant $M$.
Therefore, any conformal map 
$\varphi:\D\to D_{\eps}$ is H\"older continuous with some universal exponent $\alpha=\alpha(M)>0$;
see \cite{P}, Chapter 5.2, or \cite{GM}, Chapter 7.
Thus (\ref{e.ho})  implies 
$$|f\circ\varphi(z)-f\circ\varphi(z')| \leq C\, |z-z'|^{\alpha p/2}$$
for all $z,z'\,\in\D,$ with a constant $C$ depending on $\eps$, but with exponent 
$\alpha p/2$ depending on $\kappa$ only. 

If $\alpha\, p/2$ happened to be greater than 1/2, then the result would follow from 
the trivial estimate for the change of Hausdorff dimension by the reciprocal of the H\"older 
exponent. 
Since 
it is not the case (recall $p=8/\kappa -1$), we need to use some more advanced results.
Here is a Theorem by Jones and Makarov (Theorem C.2 in \cite{JM},
see Corollary 3.2 in \cite{KR} for a different proof of this Theorem).
\begin{thm}[Jones-Makarov]
Let $\eta\in (0,1)$ and let $\Omega$ be some H\"older domain with exponent $\eta$ (i.e., $\Omega$ is a Jordan domain
so that any conformal mapping $f:\D \rightarrow \Omega$ is $\eta$-H\"older in the disk). Then 
\begin{equation}
\dim\, \p \Omega\, \leq 2 - c\,\eta\,, \nonumber
\end{equation}
where $c>0$ is an absolute constant.
\end{thm}

Applied to our setting, this theorem implies that 
$$\dim \partial f(D_{\eps}) = \dim \partial\bigl( f\circ\varphi(\D)\bigr) \leq 2- c\alpha p/2.$$
Above, we have seen that for every $0<t<T<\infty$ we almost surely have
$\gamma[t,T] \subset D_{\eps}$ for some $\eps>0$, and therefore
$$f\bigl(\gamma(0,\infty)\cap \partial\D\bigr) \subset \bigcup_{n} \partial f(D_{1/n}).$$
Setting $B=\partial G\cap\bigcup_{n} \partial f(D_{1/n})$,
Theorem~\ref{dimension} follows at once with
$d(\kappa) = 2- c\alpha p/2$. Notice here that $B$ is indeed a $F_\sigma$ set
since by (\ref{e.ho}), $f$ is uniformly continuous on $D_{1/n}$.
\qed\medskip

 Notice that $\alpha$ can be chosen independently of $\kappa$, therefore there is 
a constant $C_1>0$ such that 
\begin{equation}\label{e.bound8}
d(\kappa)\leq 2- C_1\,(8-\kappa)\,,
\end{equation}
which proves part of Proposition \ref{p.estimate}.

\medskip

In order to show that we cannot choose $d(\kappa)\leq1$ in general, fix $\kappa\in (4,8)$
and choose $\alpha\in (\frac8{\kappa}-1,1).$ Then there is a simply connected domain $G$ whose
boundary has dimension greater than 1 and such that any conformal map $f$ from $\D$ onto $G$ is
H\"older continuous with exponent $\alpha.$ In fact, any sufficiently \lq\lq flat\rq\rq\ snowflake curve, or
the bounded Fatou component of the quadratic polynomial $z^2 +\lambda z$ with 
$0<|\lambda|<(1-\sqrt{\alpha})/(1+\sqrt{\alpha})$ will do, see \cite{AIM}, Chapter 13.3.
If $B\subset \partial G$ is a Borel set which almost surely contains
the intersection $\gamma(0,\infty)\cap\partial G$ of
the chordal  $\mathrm{SLE}_\kappa$ trace $\gamma$ from $a$ to $b$ with the boundary of the domain,
then almost surely 
the chordal $\mathrm{SLE}_\kappa$ in $\D$ does not intersect $L=\p \D \setminus f^{-1}(B)$. 
We claim that $L$ has to be of Hausdorff dimension at most $a:=\frac8{\kappa} -1$.
We briefly sketch the proof, which is based on estimates and arguments from~\cite{SZ};
see~\cite{Z} for further details.
Consider the chordal $\mathrm{SLE}_\kappa$ path from $0$ to $\infty$ in the
upper half plane. Let $L'\subset[1,2]$ be Borel-measurable and have Hausdorff dimension larger than $a$.
By Theorem 8.8 in \cite{Mattila}, there exists a Frostman measure $\mu$
supported on a compact subset $A$ of $L'$ whose $a$-energy is finite; namely, $\mu$ is supported on $A \subset L'$,
$\mu(A)>0$, and
$$\iint_{A\times A} \frac{d\mu(x) d\mu(y)} {|x-y|^a} <\infty \,.
$$
(See~\cite[Section 8]{Mattila} for background on Frostman measures.)
Let $C_\eps$ be defined as in~\cite[Section 2]{SZ}.
We now apply a second moment argument to the random variable $\mu(C_\eps)$.
By Propositions~2.3 and~2.4 in~\cite{SZ} we have for $1\le x< y\le 2$ and $\eps<1$ that
$\P\bigl[x\in C_\eps\bigr]$ is comparable to  $\eps^a$ and 
$\P\bigl[x,y\in C_\eps\bigr]\le C\,\eps^{2a}\,(y-x)^{-a}$.
Hence,
$$
\E\bigl[\mu(C_\eps)\bigr]=\int_{[1,2]} \P\bigl[x\in C_\eps\bigr]\,d\mu(x)
$$
is of order $\eps^a$ and
\begin{multline*}
\E\bigl[\mu(C_\eps)^2\bigr] = \iint_{[1,2]^2} \P\bigl[x,y\in C_\eps]\,d\mu(x)\,d\mu(y)\\
\le C\,\eps^{2a}\,\iint_{[1,2]^2} |x-y|^{-a}\,d\mu(x)\,d\mu(y)\,.
\end{multline*}
By the choice of $\mu$, the latter is bounded by a constant times $\eps^{2a}$. Thus,
$\E\bigl[\mu(C_\eps)^2\bigr] \le C\, \E\bigl[\mu(C_\eps)\bigr]^2$.
The standard second moment argument (i.e., Cauchy Schwarz) therefore implies
that $\P[\mu(C_\eps)>0]$ is bounded away from $0$ independently of $\eps$.
Thus, with positive probability $\mu(C_\eps)>0$ for every $\eps>0$.
But since the support of $\mu$ is compact and contained in $L'$, this implies that
the SLE path hits $L'$ with positive probability, which clearly implies our claim
that the Hausdorff dimension of $L=\partial \D\setminus f^{-1}(B)$ is at most $a$.

As $\dim L \le a$ and $f$ is $\alpha$-H\"older, it follows that $\dim f(L)\le a/\alpha<1$.
Since $B=\partial G\setminus f(L)$, it follows that $\dim B=\dim \partial G>1$, as required.

\bigskip

\noindent{\bf Proof of Proposition \ref{p.estimate}.}
We already noticed one part of the inequality in~\eqref{e.bound8}.
It remains to bound $d(\kappa)$ when $\kappa$ is close to 4. We will follow the same plan of proof, but instead of using
the quantities $\int_{Q} |f'|^2$ we will refine Lemma \ref{content} by using $\int_{Q} |f'|^t$ for some well chosen $t=t(\kappa)$
close to zero. What allowed us to conclude the proof of Lemma \ref{content} was the fact that for bounded domains, $\int_{\D} |f'|^2 <\infty$.
Here we will use instead the known bounds on the Integral means spectrum of univalent functions
(in particular, we do not need to assume $G$ bounded in this proof).

Let us briefly recall some facts about the {\it Integral Means Spectrum} (see \cite{P}).
Let $f$ be an univalent function in the unit disc $\D$. For any $t\in \R$, let
$$
\beta_f(t):=\inf\Bigl\{\beta\in\R: \lim_{r\to 1} (1-r)^\beta\int_{|z|=r} |f'(z)|^t\, |dz|=0
\Bigr\}.
$$
The {\it universal integral means spectrum} $B(t)$ of univalent functions is defined as
\begin{equation}
B(t) = \sup_f \beta_f (t),\nonumber
\end{equation}
where the supremum is over all univalent functions (often one restricts the supremum to bounded
univalent functions, resulting in a slightly different spectrum). 
Much is known about this spectrum, see \cite{P} and references therein.
We will use the following upper bound on the spectrum (\cite{P}, Theorem 8.5):
\begin{equation}
B(t) < 4 t^2 \, \text{ for all }t\in \R\setminus\{0\}\,. \nonumber
\end{equation}

We prove the following Lemma (which is a refinement of Lemma \ref{content} for $p$ close to 1), from which 
Proposition \ref{p.estimate} will easily follow. 
\begin{lemma}\label{content2}
Let $G\subsetneqq \C$ be a simply connected domain and $f\,:\,\D\rightarrow  G$ a 
conformal homeomorphism. For every $0<p< 1$, and $\eps>0$ there is $D \subset \D$
and $C>0$ such that 
\begin{equation}
H_p(\D \setminus D) <\eps, \nonumber
\end{equation}
and
\begin{equation}\label{e.holder2}
|f(z)-f(z')| \leq C |z - z'| ^{1-6\sqrt{1-p}},
\end{equation}
for all $z,z'\in D$.
\end{lemma}

\proof
Notice that the lemma is relevant only when $p$ is close enough to 1.
As in the proof of Lemma \ref{content}, we consider the collection $\mathcal{Q}$ of dyadic squares.
For each $Q\in \mathcal{Q}$, denote by $L(Q)$ the inner ``segment'' $\{z\in Q\,:\, |z|= 1-2^{-n} \}$.
Fix the parameters $N>1, t >0$ and $0<\delta< p$, to be determined later, and let $\mathcal{L}$ be the collection of squares 
$Q\in \mathcal{Q}$ for which $$l(Q) \leq 2^{-N}$$ and $$\int_{L(Q)} |f'|^{t} > l(Q)^{p-\delta}.$$
Set as before $$ D= \D \setminus \bigcup_{Q\in \mathcal{L}} Q\,.$$
Using the integral means spectrum, one can easily estimate the Hausdorff $p$-content of $\D\setminus D$. Indeed,

\begin{eqnarray}
H_p(\D\setminus D) \,\leq\,  \sum_{Q\in \mathcal{L}} l(Q)^p
&\leq & \sum_{Q\in \mathcal{L}} l(Q)^{\delta} \int_{L(Q)} |f'|^t \nonumber\\
&\leq& \sum_{n\geq N } 2^{-n\delta} \int_{|z|= 1-2^{-n}} |f'|^t \nonumber\\
&\leq & \sum_{n\geq N} C\, 2^{-n(\delta - 4t^2)},\nonumber
\end{eqnarray}
since $\beta_f(t)\leq B(t) < 4t^2$, where $C>0$ may depend on $f$.
Therefore one needs to choose $\delta> 4t^2$; let $\delta:= 5t^2$. By taking $N$ large enough,
we get $H_p(\D\setminus D) \leq \eps$.

As in Lemma \ref{content}, to conclude the proof, it is enough to check that there is some $C=C(f,N,p)$
such that for all $z\in D$,
\begin{equation}\label{e.claim2}
|f'(z)| \leq C\, \frac 1 {(1-|z|)^{6 \sqrt{1-p}}}\,.
\end{equation}
Let $z\in D,\, |z|\geq 1-2^{-N}$ (for $|z|< 1-2^{-N}$, we choose $C$ large enough so that (\ref{e.claim2}) is satisfied).
Let $Q$ be such that $z\in T(Q)$ and notice that $Q\notin \mathcal{L}$.
Therefore $\int_{L(Q)} |f'|^t \leq l(Q)^{p-5t^2}$. By Koebe's distortion
theorem, $|f'|$ fluctuates by at most a multiplicative constant with is essentially constant 
within $T(Q)$, and hence
$$
|f'(z)| \leq O(1)\, \Bigl(\frac {1}{1-|z|}\Bigr)^{\frac{1-p + 5t^2}{t}}\,.
$$
 By choosing our last parameter $t:= \sqrt{1-p}$,
this leads to (\ref{e.claim2}).

Now, following Lemma \ref{content}, i.e., integrating along appropriate arcs, this proves Lemma \ref{content2}.
\qed
\medskip

Proposition \ref{p.estimate} follows from Lemma \ref{content2} with $p=8/\kappa -1$ in the following way: a.s.\ 
the $\mathrm{SLE}_{\kappa}$ trace remains in $D=D_{\eps}$ for some $\eps$ small enough. Moreover the map $f$ from
$D_{\eps}$ to $f(D_{\eps})$ is $\eta$-H\"older with $\eta = 1- 6\sqrt{1-p}$.
Hence, by the obvious bound (here we do not need the above Theorem of Jones and Makarov),
\begin{equation}
\dim \partial f(D_{\eps}) \leq  \frac 1 {\eta} \dim \partial D_{\eps} = \frac 1 {1 - 6\sqrt{1-p}} \leq 1 +  C_2 \sqrt{\kappa - 4}\nonumber
\end{equation}
($p=8/\kappa -1$), which together with~\eqref{e.bound8} implies Proposition \ref{p.estimate}.
\qed

\section{Related results}

In \cite{SZ,AS},
it is proved that for $\kappa\in(4,8]$,
the chordal $\mathrm{SLE}_{\kappa}$ in $\H$ a.s.\ satisfies $\dim(\mathrm{SLE}_\kappa\cap \R)= 2-8/\kappa$
(the same holds of course for radial $\mathrm{SLE}$ in the unit disc $\D$). What happens in the case of a general 
simply connected domain $G$? 
That is, what is the Hausdorff dimension of $\mathrm{SLE}_{\kappa}\cap \partial G$ ?
By Theorem \ref{dimension} we know that it is $\leq d(\kappa)<2.$
If $G$ is a John domain, we can generally do much better.
Notice that if one takes $\kappa=8$, 
then $\mathrm{SLE}_{\kappa}$ is space filling and thus
$\mathrm{SLE}_{\kappa}\cap \partial G =\partial G.$ 
If $G$ is a John domain, it is known (\cite{P}, Theorem 10.17)
that the dimension $d$ of the boundary
is the unique solution to the equation $\beta_f(d)= d-1,$
where $f$ is a conformal map from $\D$ to $G$
and $\beta_f(d)$ is the integral means spectrum of $f$.
We will now sketch a proof that for general $\kappa\in(4,8]$
and John domains $G$,
the dimension of $\mathrm{SLE}_{\kappa}\cap \partial G$ 
is bounded from above by the solution $d$ of the equation $\beta_f(d)=d-(2-8/\kappa)$.
Fix $\kappa\in(4,8)$,
and $t>0$ such that $\beta_f(t)<t-(2-8/\kappa)$.
Let $\gamma$ be a chordal $\mathrm{SLE}_{\kappa}$ from 
$-1$ to 1 in $\D$. We are interested in $\dim f(\gamma(0,\infty))\cap \partial G$. 
Since we assumed $G$ to be a John domain, there is a constant $C$ such that 
$C\cdot f(T(Q))\supset f(Q)$
for each dyadic square $Q\in\mathcal{Q}$.
Cover $f(\gamma(0,\infty))\cap \partial G$ by $\bigcup C\cdot f(T(Q)),$
where the union is over those $Q\in\mathcal{Q}$ 
for which $\gamma\cap 2\cdot Q \neq \emptyset$ 
and $l(Q)\leq 2^{-N}$ with $N$ large enough.
By~\eqref{hitting}, the expected Hausdorff $t$-content of $f(\gamma(0,\infty))\cap \partial G$ is thus bounded by
\begin{equation}
\sum_{n\geq N} \sum_{l(Q)=2^{-n}} O(1) (2^{-n} |f'(z_Q)|)^t (2^{-n})^{8/\kappa -1} = \sum_{n\geq N}O(1) 2^{-n(t-2+8/\kappa)} \int_{|z|=1-2^{-n}} |f'|^t,\nonumber
\end{equation}
where $z_Q$ is any point in $T(Q).$ 
This sum converges since $\beta_f(t)< t-(2-8/\kappa)$; so by letting $N$ going to $\infty$, the expected Hausdorff $t$-content of $f(\gamma(0,\infty))\cap \partial G$
is equal to zero, and the result follows. 

%%%%%%%%%%%%%%%%%%
\begin{figure}[htbp]
\centerline{\epsfysize=3 true in \epsffile{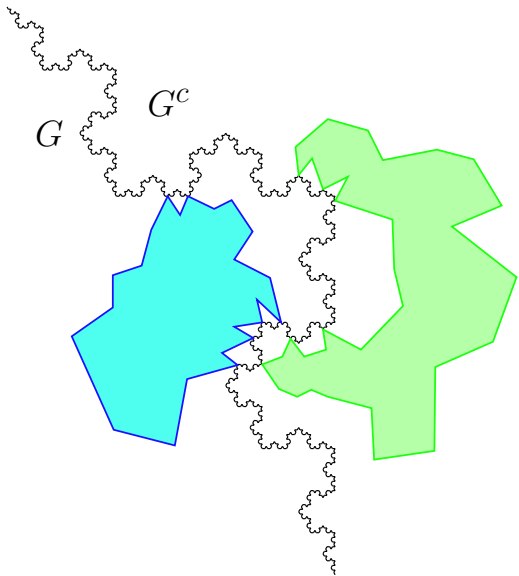}} 
\caption{\label{f.singularity} For some domains $G$, (at the continuum limit) percolation clusters inside $G$ are ``invisible'' to percolation clusters 
inside $G^c$.}
\end{figure}
%%%%%%%%%%%%%%%%

\bigskip
If $\kappa > 16/3,$ then $\dim(\mathrm{SLE}_\kappa\cap \R)= 2-8/\kappa > 1/2$ and hence two
independent $\mathrm{SLE}_\kappa,$ one in the upper half plane and one in the lower, will intersect a.s.
This is not true any more for general Jordan domains: For each 
$\kappa\in (4,8)$, there exists a John domain (actually a quasidisc) $G=G(\kappa)$ and a set 
$E\subset \p G$ such that if $\gamma_1$ and $\gamma_2$ are respectively  $\mathrm{SLE}_{\kappa}$ curves driven inside and outside $G$,
then a.s.\ $\gamma_1(0,\infty)\cap \p G \subset E$ while $\gamma_2(0,\infty)\cap \p G \subset E^c$.
Indeed, given $\kappa\in (4,8)$ and choosing $0<\varepsilon <8/\kappa-1$,
by \cite{T} and \cite{R} there is a quasidisc $G$ and a subset $A\subset \partial\D$ with
$\dim A<\varepsilon$ and $\dim \partial\D \setminus f_c^{-1}(f(A))<\varepsilon,$
where $f$ and $f_c$ are conformal maps from $\D$ to $G$ and $G^c.$ It follows that 
$\mathrm{SLE}_\kappa$ in $\D$ will a.s. be disjoint from both $A$ and 
$\partial\D \setminus f_c^{-1}(f(A))$, and the claim follows with $E=f(\partial\D \setminus A)$. 
This can be viewed as an $\mathrm{SLE}$ analog of the Theorem by 
Bishop, Carleson, Garnett and Jones about harmonic measure (see \cite{BCGJ, R}). Figure \ref{f.singularity}
is an illustration of this property in the case of percolation clusters.


\begin{thebibliography}{abc}

\bibitem{AK} T. Alberts, M. Kozdron,
\newblock Intersection probabilities for a chordal SLE path and a semicircle.
\newblock {\it Electron. Comm. Probab.} (2008), 13:448--460. 


\bibitem{AS} T. Alberts, S. Sheffield,
\newblock Hausdorff dimension of the SLE curve intersected with the real line.
\newblock {\it Electron. J. Probab.} (2008), 13:1166--1188.


\bibitem{AIM} K. Astala, T. Iwaniec, G. Martin,
Elliptic Partial Differential Equations and
Quasiconformal Mappings in the Plane, to appear.

\bibitem{BCR} I. Benjamini, Z. Chen, S. Rohde,
Boundary Trace of Reflecting Brownian Motions.
{\it Probab. Theory Relat. Fields \bf 129} (2004), 1--17.

\bibitem{BCGJ}
C. Bishop, L. Carleson, J. Garnett and P. Jones,
\newblock Harmonic measures supported on curves.
\newblock {\it Pacific J. Math. \bf 138} no. 2 (1989), 233--236.


\bibitem{GM} J. Garnett, D. Marshall, Harmonic measure.
Cambridge University Press (2005)

\bibitem{JM} P. Jones, N. Makarov, Density properties of harmonic measure.
{\it Ann. Math. \bf 142} (1995), 427-455.

\bibitem{Ke} H. Kesten, The incipient infinite cluster in two-dimensional percolation.
{\it Probab. Th. Rel. Fields \bf 73} (1986), 369-394.

\bibitem{KR} P. Koskela, S. Rohde, Hausdorff dimension and mean porosity. 
{\it Math. Annalen \bf 309} (1997), 593--609. 

\bibitem{K} J. Kuester, A region whose prime ends all have the same impression.
{\it Math. Z. \bf 136} (1974), 1--5.

\bibitem{Mattila} P. Mattila, Geometry of sets and measures in Euclidean spaces. Cambridge University Press 1995.

\bibitem{P} C. Pommerenke, Boundary behaviour of conformal maps. Springer
Verlag, Berlin Heidelberg 1992.

\bibitem{R} S. Rohde, On conformal welding and quasicircles. 
{\it Michigan Math. J. \bf 38} (1991), 111--116.

\bibitem{RS} S. Rohde, O. Schramm, Basic properties of SLE.
{\it Ann. Math. \bf 161} (2005), 879--920.

\bibitem{S} O. Schramm, Conformally invariant scaling limits:
an overview and a collection of problems. In {\it Proceedings of ICM
Madrid 2006, \bf  1}, 513--543. European Math. Soc. 2007.

\bibitem{SheffCLE} Scott Sheffield,
Exploration trees and conformal loop ensembles. To appear in
{\it Duke Mathematical Journal}, arXiv:math.PR/0609167.

\bibitem{SZ} O. Schramm, W. Zhou, Boundary proximity of SLE. arXiv:0711.3350v3

\bibitem{T} P. Tukia, Hausdorff dimension and quasisymmetric mappings.
{\it Math. Scand \bf 65} (1989), 152--160.

\bibitem{Z} W. Zhou, in preparation.

\end{thebibliography}
\end{document}